\documentclass[12pt,a4paper]{article}
\usepackage[latin 2]{inputenc}
\usepackage{amssymb}
\usepackage{latexsym}
\usepackage{amsrefs}
\usepackage{amsmath}
\usepackage{amsthm}
\newtheorem{thm}{Theorem}
\newtheorem{prop}[thm]{Proposition}
\newtheorem{lemma}[thm]{Lemma}
\newtheorem{definition}{Definition}\theoremstyle{definition}
\begin{document}

\begin{center}
{\large\textbf{Lattice points in algebraic cross-polytopes and simplices}}

\vspace{5mm}

\textbf{Bence Borda}

{\footnotesize Department of Mathematics, Rutgers University

110 Frelinghuysen Road, Piscataway, NJ-08854, USA

Email: \texttt{bordabence85@gmail.com}}

\vspace{5mm}


\end{center}

\vspace{5mm}

\begin{abstract}

The number of lattice points $\left| tP \cap \mathbb{Z}^d \right|$, as a function of the real variable $t>1$ is studied, where $P \subset \mathbb{R}^d$ belongs to a special class of algebraic cross-polytopes and simplices. It is shown that the number of lattice points can be approximated by an explicitly given polynomial of $t$ depending only on $P$. The error term is related to a simultaneous Diophantine approximation problem for algebraic numbers, as in Schmidt's theorem. The main ingredients of the proof are a Poisson summation formula for general algebraic polytopes, and a representation of the Fourier transform of the characteristic function of an arbitrary simplex in the form of a complex line integral.

\end{abstract}

\vspace{5mm}

\section*{Acknowledgment}

This paper is based on the doctoral dissertation of the author. The author is grateful to his advisor, J\'ozsef Beck.

\section{Introduction}

Given a set $P \subset \mathbb{R}^d$, estimating the number $\left| tP \cap \mathbb{Z}^d \right|$ of lattice points in its dilates
\[ tP = \left\{ tx \,\, \middle| \,\, x \in P \right\} , \]
as a function of the real variable $t>1$ is a classical problem in number theory. The case when $P$ is a convex body with a smooth boundary has a vast literature, and will not be considered in this paper. Instead, we shall study the case when $P$ is a polytope, i.e. the convex hull of finitely many points in $\mathbb{R}^d$. Moreover, we shall focus on polytopes $P$ defined in terms of algebraic numbers.

There is an important class of such polytopes for which the lattice point counting problem is completely solved. If every vertex of the polytope $P \subset \mathbb{R}^d$ is a lattice point, and $P$ has a nonempty interior, then there exists a polynomial $p(t) \in \mathbb{Q}[t]$ of degree $d$ such that
\[ \left| tP \cap \mathbb{Z}^d \right| = p(t) \]
for every positive integer $t$. This is \textbf{Ehrhart's theorem} \cites{Ehr1,Ehr2,Ehr3}, and the polynomial $p(t)$ is called the Ehrhart polynomial of $P$. It is also known that the leading coefficient of $p(t)$ is the Lebesgue measure of $P$, while the coefficient of $t^{d-1}$ is one half of the normalized surface area of the boundary $\partial P$. Here the normalized surface area of a $d-1$ dimensional face of $P$ is defined as the surface area of the face divided by the covolume of the $d-1$ dimensional sublattice of $\mathbb{Z}^d$ on the affine hyperplane containing the face.

Ehrhart's theorem can actually be generalized to polytopes with vertices in $\mathbb{Q}^d$ instead of $\mathbb{Z}^d$. Moreover, we can allow the dilation factor $t$ to be a positive rational or real number. In this more general case there still exists a precise formula without any error term for the number of lattice points in $tP$, in the form of a so-called quasi-polynomial \cites{Bald, Linke}. Not surprisingly, the coefficients of these Ehrhart quasi-polynomials depend on the fractional part of certain integral multiples of $t$.

There is no complete answer to the lattice point counting problem, however, if we only assume that the vertices of the polytope $P$ have algebraic coordinates. The first result regarding this more general case is due to Hardy and Littlewood \cites{Hardy1,Hardy2}. Let
\begin{equation}\label{triangle}
S = \left\{ (x,y) \in \mathbb{R}^2 \,\, \middle| \,\, x,y \ge 0, \,\, \frac{x}{a_1} + \frac{y}{a_2} \le 1 \right\} ,
\end{equation}
i.e. the closed right triangle with vertices $(0,0) , (a_1 , 0) , (0, a_2)$, where $a_1 , a_2 >0$. As observed by Hardy and Littlewood, estimating $\left| tS \cap \mathbb{Z}^d \right|$ for real numbers $t>1$ is closely related to the classical Diophantine problem of approximating the slope $- \frac{a_2}{a_1}$ by rational numbers with small denominators. If the slope $- \frac{a_2}{a_1}$ is algebraic, then
\begin{equation}\label{eq1}
\left| tS \cap \mathbb{Z}^2 \right| = \frac{a_1 a_2}{2} t^2 + \frac{a_1 + a_2}{2} t + O \left( t^{\beta} \right)
\end{equation}
for some $0< \beta <1$ depending only on $a_1 , a_2$. This groundbreaking theorem was one of the first results on Diophantine approximation of general algebraic numbers. Note that the main term in \eqref{eq1} is a polynomial, where the leading coefficient is the area of $S$, while the coefficient of $t$ is one half of the total length of the legs of the right triangle $S$.

Later Skriganov \cite{Skrig2} studied the lattice point counting problem in more general polygons whose sides have algebraic slopes. From his results it follows easily that the error term in \eqref{eq1} can be improved to $O \left( t^{\varepsilon} \right)$ for any $\varepsilon >0$. His main idea was to combine the Poisson summation formula and \textbf{Roth's theorem}
\[ \inf_{m>0} m^{1+\varepsilon} \left\| m \alpha \right\| > 0, \]
applied to the algebraic slopes of the sides of the polygon. Note that throughout the paper $|\cdot |$ denotes the Euclidean norm of a real number or vector, or the cardinality of a set, while $\left\| \cdot \right\|$ is the distance from the nearest integer function.

In the special case when the slope $-\frac{a_2}{a_1}$ is a quadratic irrational, then \eqref{eq1} in fact holds with an error term $O \left( \log t \right)$, which is actually best possible. This observation was already made by Hardy and Littlewood \cites{Hardy1,Hardy2}, and is related to the fact that the Diophantine approximation problem for quadratic irrationals is much easier, than it is for general algebraic numbers.

Much less is known about higher dimensional lattice point counting problems. Trivially, for any polytope $P \subset \mathbb{R}^d$ we have
\begin{equation}\label{eq2}
\left| tP \cap \mathbb{Z}^d \right| = \lambda (P) t^d + O \left( t^{d-1} \right)
\end{equation}
with an implied constant depending only on $P$, where $\lambda (P)$ denotes the Lebesgue measure of $P$. In a sense \eqref{eq2} is best possible. Indeed, consider the normal vectors of the $d-1$ dimensional faces of $P$. Here and from now on by a normal vector of a $d-1$ dimensional face we mean any nonzero vector orthogonal to the face, not necessarily of unit length. It is easy to see that if $P$ contains the origin in its interior, and it has a $d-1$ dimensional face with a rational normal vector, then
\[ \left| tP \cap \mathbb{Z}^d \right| = \lambda (P) t^d + \Omega \left( t^{d-1} \right) . \]

Partial results have been obtained in the case when the polytope $P$ is subjected to certain irrationality conditions. \textbf{Randol's theorem} \cite{Rand} states that if every $d-1$ dimensional face of a polytope $P \subset \mathbb{R}^d$ has a normal vector with two coordinates of algebraic irrational ratio, then \eqref{eq2} holds with an error term $O \left( t^{d-2+ \varepsilon} \right)$ for any $\varepsilon >0$. The proof is again based on the Poisson summation formula and Roth's theorem applied to the algebraic ratios.

Skriganov \cite{Skrig1} introduced methods of ergodic theory in lattice point counting problems with respect to more general lattices. For certain pairs of algebraic polytopes $P$ and algebraic unimodular lattices $\Gamma$ it is proved \cite{Skrig1}*{Theorem 2.3} that
\[ \left| t P \cap \Gamma \right| = \lambda (P) t^d + O \left( t^{\varepsilon} \right) \]
for any $\varepsilon >0$.

Stronger results have been obtained in the case when a random translation and/or random rotation, in the sense of the Haar measure on $\textrm{SO}(d)$, is applied to a polytope \cites{Brand,Skrig1,Tarn}. Since a randomly translated or rotated polytope loses any kind of algebraicity, these results are outside the scope of this paper.

\section{Main results}

\subsection{Statement of the problems}\label{section2.1}

In the present paper we wish to study the lattice point counting problem in two specific polytopes. Let $d \ge 2$, $a_1 , \dots , a_d >0$, and consider
\begin{equation}\label{cross}
C = C (a_1 , \dots , a_d) = \left\{ x \in \mathbb{R}^d \,\, \middle| \,\, \frac{|x_1|}{a_1} + \cdots + \frac{|x_d|}{a_d} \le 1 \right\} ,
\end{equation}
\begin{equation}\label{simplex}
S = S(a_1 , \dots , a_d) = \left\{ x \in \mathbb{R}^d \,\, \middle| \,\, x_1 , \dots , x_d \ge 0, \,\, \frac{x_1}{a_1} + \cdots + \frac{x_d}{a_d} \le 1 \right\} .
\end{equation}

Here $C$ is a cross-polytope whose $d-1$ dimensional faces have normal vectors of the form
\[ \left( \frac{\pm 1}{a_1} , \dots , \frac{\pm 1}{a_d} \right) . \]
The vertices of $C$, on the other hand, are of the very simple form $\left( 0 , \dots , \pm a_i , \dots , 0 \right)$ for some $1 \le i \le d$. The polytope $S$ is a simplex the vertices of which are the origin and the points $\left( 0 , \dots , a_i , \dots , 0 \right)$ for $1 \le i \le d$. Note that $S$ is a direct generalization of the right triangle \eqref{triangle} studied by Hardy and Littlewood.

We wish to study $\left| tC \cap \mathbb{Z}^d \right|$ and $\left| tS \cap \mathbb{Z}^d \right|$, as $t \to \infty$ along the reals under the assumption that $\frac{1}{a_1} , \dots , \frac{1}{a_d}$ are algebraic and linearly independent over $\mathbb{Q}$. Our main result is that there exist explicitly computable polynomials $p(t)$ and $q(t)$ such that
\[ \begin{split} \left| tC \cap \mathbb{Z}^d \right| &= p(t) + O \left( t^{\frac{(d-1)(d-2)}{2d-3} + \varepsilon} \right) , \\ \left| tS \cap \mathbb{Z}^d \right| &= q(t) + O \left( t^{\frac{(d-1)(d-2)}{2d-3} + \varepsilon} \right) \end{split} \]
for any $\varepsilon >0$. For the precise formulation of the main results see Theorems \ref{theorem6}, \ref{theorem7} and \ref{theorem8} in Section \ref{section2.4}.

We start with the simple observation that these two problems are equivalent.

\begin{prop}\label{proposition1}
Let $a_1 , \dots , a_d >0$ be arbitrary reals, and let $S$ be as in \eqref{simplex}. For every $I \subseteq [d] = \left\{ 1 , 2 , \dots , d \right\}$ let
\[ C_I = \left\{ x \in \mathbb{R}^d \,\, \middle| \,\, \sum_{i \in I} \frac{|x_i|}{a_i} \le 1 , \,\, \forall j \in [d] \backslash I : \,\, x_j=0 \right\}. \]
Then for any real $t>0$ we have
\[ \left| tS \cap \mathbb{Z}^d \right| = \frac{1}{2^d} \sum_{I \subseteq [d]} \left| t C_I \cap \mathbb{Z}^d \right| . \]
\end{prop}

\noindent\textbf{Proof:} For every $\sigma \in \left\{ 1 , -1 \right\}^d$ consider the simplex
\begin{equation}\label{triangulation}
S_{\sigma} = \left\{ x \in \mathbb{R}^d \,\, \middle| \,\, \sigma_1 x_1 \ge 0 , \dots , \sigma_d x_d \ge 0 , \,\, \frac{\sigma_1 x_1}{a_1} + \cdots + \frac{\sigma_d x_d}{a_d} \le 1 \right\} .
\end{equation}
We have
\begin{equation}\label{sum}
\sum_{\sigma \in \left\{1, -1 \right\}^d} \left| tS_{\sigma} \cap \mathbb{Z}^d \right| = \sum_{I \subseteq [d]} \left| t C_I \cap \mathbb{Z}^d \right| .
\end{equation}
Indeed, a lattice point in $tC \cap \mathbb{Z}^d$ with $k$ zero coordinates is counted $2^k$ times on both sides of \eqref{sum}. Finally, note that the sum on the left hand side of \eqref{sum} has $2^d$ terms, and that each term equals $\left| tS \cap \mathbb{Z}^d \right|$.
\begin{flushright}
$\square$
\end{flushright}

It should be noted that Skriganov \cite{Skrig1}*{Theorem 6.1} proved a quite general bound for the lattice discrepancy
\[ \left| tP \cap \mathbb{Z}^d \right| - \lambda (P) t^d \]
for an explicitly defined, wide class of polytopes, which in a sense contains ``almost every'' polytope. One can check, however, that neither $C$, nor $S$ belongs to this wide class.

The rest of the paper is organized as follows. In Section \ref{section2.2} we introduce a Poisson summation formula for algebraic polytopes. A new representation of the Fourier transform of the characteristic function of an arbitrary simplex in $\mathbb{R}^d$ is given in Section \ref{section2.3}. The main results of the paper are stated in Section \ref{section2.4}, while conclusions are listed in Section \ref{section2.5}. Finally, the proofs of all the results are given in Section \ref{section3}.

\subsection{Poisson summation formula for algebraic polytopes}\label{section2.2}

Given a polytope $P \subset \mathbb{R}^d$ and a real number $t>0$, let $\chi_{tP}$ denote the characteristic function of $tP$, and let
\[ \hat{\chi}_{tP} (y) = \int_{tP} e^{-2 \pi i \langle x,y \rangle} \, \mathrm{d}x \]
denote its Fourier transform, where $\langle x,y \rangle$ is the scalar product of $x,y \in \mathbb{R}^d$. The main idea is to apply the Poisson summation formula
\begin{equation}\label{Poisson}
\left| tP \cap \mathbb{Z}^d \right| = \sum_{m \in \mathbb{Z}^d} \chi_{tP} (m) \sim \sum_{m \in \mathbb{Z}^d} \hat{\chi}_{tP} (m).
\end{equation}
Here the symbol $\sim$ means that the series of Fourier transforms in \eqref{Poisson} has to be treated as a formal series, which may or may not converge. The reason for this is that the Poisson summation formula only holds for sufficiently smooth functions, and $\chi_{tP}$ is not even continuous. To ensure convergence we introduce the Ces\`aro means of the series as follows.

\begin{definition}\label{definition1}
For a polytope $P \subset \mathbb{R}^d$, a real number $t>0$ and an integer $N>0$ let
\[ \mathrm{Ces} (tP,N)= \frac{1}{N^d} \sum_{M \in [0,N-1]^d} \sum_{m \in [-M_1 , M_1] \times \cdots \times [-M_d , M_d]} \hat{\chi}_{tP} (m) . \]
\end{definition}

The number of lattice points in $tP$ can be approximated by the Ces\`aro means using the following theorem.

\begin{thm}[Poisson summation formula for algebraic polytopes]\label{theorem2}
Let $P \subset \mathbb{R}^d$ be a polytope with a nonempty interior, and let $2 \le k \le d$. Suppose that every $d-1$ dimensional face of $P$ has a normal vector $(n_1 , \dots , n_d)$ such that its coordinates are algebraic and span a vector space of dimension at least $k$ over $\mathbb{Q}$. Then for every real $t>1$, every integer $N>1$ and every $\varepsilon >0$ we have
\[ \left| tP \cap \mathbb{Z}^d \right| = \mathrm{Ces} (tP,N) + O \left( t^{d-k} + t^{d-1+\varepsilon} \sqrt{\frac{\log N}{N}} \right) . \]
The implied constant depends only on $P$ and $\varepsilon$, and is ineffective.
\end{thm}

Note that under the assumptions of Theorem \ref{theorem2} it is possible that the affine hyperplane containing a $d-1$ dimensional face of $tP$ contains a $d-k$ dimensional sublattice of $\mathbb{Z}^d$, as $t \to \infty$ along a special sequence. Thus if we are to approximate $|tP \cap \mathbb{Z}^d |$ by any continuous function, an error of $t^{d-k}$ is inevitable. This inevitable error is minimized by assuming $k=d$, i.e. that the coordinates of the normal vectors are algebraic and linearly independent over $\mathbb{Q}$.

The proof of Theorem \ref{theorem2} is based on \textbf{Schmidt's theorem} \cite{Schmidt}, which states that if $\alpha_1 , \dots , \alpha_d$ are algebraic reals such that $1 , \alpha_1 , \dots , \alpha_d$ are linearly independent over $\mathbb{Q}$, then
\begin{equation}\label{Schmidt1}
\inf_{m \in \mathbb{Z}^d \backslash \{ 0 \}} |m|^{d+ \varepsilon} \left\| m_1 \alpha_1 + \cdots + m_d \alpha_d \right\| >0 ,
\end{equation}
and
\begin{equation}\label{Schmidt2}
\inf_{m>0} m^{1 + \varepsilon} \left\| m \alpha_1 \right\| \cdots \left\| m \alpha_d \right\| >0
\end{equation}
for any $\varepsilon >0$. It is worth noting that we shall apply \eqref{Schmidt1} to $k-1$ algebraic numbers, where $k$ is as in Theorem \ref{theorem2}. In fact, in the most important case $k=d$ we shall apply \eqref{Schmidt1} to $\alpha_1 = \frac{n_1}{n_d} , \dots , \alpha_{d-1} = \frac{n_{d-1}}{n_d}$, and other similar pairwise ratios of the coordinates of a normal vector. The ineffectiveness of Theorem \ref{theorem2} is of course caused by the ineffectiveness of Schmidt's theorem.

It should be mentioned that in lattice point counting problems convergence in the Poisson summation formula is traditionally ensured by convolving the characteristic function by a smooth approximate identity $\eta$ with a compact support. Such a convolution only changes the values of $\chi_{tP}$ close to the boundary of $tP$, the cutoff distance being the diameter $h$ of the support of $\eta$. The error of replacing $\chi_{tP}$ by the convolution in the left hand side of \eqref{Poisson} is therefore bounded by the number of lattice points close to the boundary of $tP$, and so it can be estimated by Lemma \ref{lemma9} below. The smoothness of $\eta$ ensures that the convolution satisfies the Poisson summation formula. Moreover, $\hat{\eta}(m)$ is close to $1$ when $|m|$ is not too large, the cutoff again being related to the diameter $h$. This way we could obtain an alternative approximation for the number of lattice points in $tP$, similar to Theorem \ref{theorem2}. The limit $N \to \infty$ in Theorem \ref{theorem2} would correspond to letting the diameter $h$ approach zero.

\subsection{The Fourier transform of the characteristic function of a polytope}\label{section2.3}

In order to use the Ces\`aro means in Definition \ref{definition1} to approximate $\left| tP \cap \mathbb{Z}^d \right|$, we need to find the Fourier transform of the characteristic function of a polytope. Several authors have found explicit formulas for the case of an arbitrary polytope using the divergence theorem (e.g. \citelist{\cite{Rand}\cite{Skrig1}*{Lemma 11.3}}). The following representation, however, is a new result.

\begin{thm}\label{theorem3}
Let $S \subset \mathbb{R}^d$ be an arbitrary simplex with vertices $v_1 , \dots , v_{d+1}$. For any real $t>0$, any $y \in \mathbb{R}^d$ and any $R> \max_{j} \left| \langle v_j , y \rangle \right|$ we have
\[ \hat{\chi}_{tS} (y) = \frac{(-1)^d d!}{(2 \pi i)^{d+1}} \lambda (S) \int_{|z|=R} \frac{e^{- 2 \pi i z t}}{\left( z - \langle v_1 ,y \rangle \right) \cdots \left( z - \langle v_{d+1} ,y \rangle \right)} \, \mathrm{d}z . \]
\end{thm}

The slightly ambiguous notation $|z|=R$ in Theorem \ref{theorem3} means a complex line integral along the positively oriented circle of radius $R$ centered at the origin. The condition $R> \max_{j} \left| \langle v_j , y \rangle \right|$ ensures that every pole of the meromorphic integrand lies inside this circle.

First of all note that finding $\hat{\chi}_{tP}$ for an arbitrary polytope $P$ can be reduced to Theorem \ref{theorem3} by triangulating $P$ into simplices. It is also worth mentioning that the variable $t$ appears only in the complex exponential function in the numerator. Thus Theorem \ref{theorem3} can be regarded as a Fourier expansion of $\hat{\chi}_{tS}(y)$ in the variable $t$, with the ``frequencies'' being the points of the circle $|z|=R$.

Why is Theorem \ref{theorem3} important, especially since explicit formulas for $\hat{\chi}_{tS}(y)$ have already been known? The main advantage is that the formula in Theorem \ref{theorem3} holds for \textit{any} $y \in \mathbb{R}^d$. To apply the Poisson summation formula, we need to sum $\hat{\chi}_{tS} (y)$ over lattice points $y=m \in \mathbb{Z}^d$. Nothing prevents the poles $\langle v_j , m \rangle$ from coinciding, in which case the integrand has a higher order pole. We will apply the residue theorem to handle such cases. Note that the residue of the integrand at a high order pole contains a high order derivative of $e^{-2 \pi i z t}$ with respect to $z$, which in turn yields a high power of $t$. We shall thus use the intuition that the residues of the high order poles of the integrand in Theorem \ref{theorem3} yield the main term in the Poisson summation formula, while the residues of the simple poles yield an error term. The most extreme case of course is that of $m=0 \in \mathbb{Z}^d$, for which the integrand has a pole of order $d+1$ with residue $\lambda (S) t^d$.

Consider now the special case of the cross-polytope $C$, as in \eqref{cross}. The simplices $S_{\sigma}$, as in \eqref{triangulation}, $\sigma \in \left\{ 1 , -1 \right\}^d$, triangulate $C$ into $2^d$ simplices to which we can apply Theorem \ref{theorem3}. Since the vertices $v_1 , \dots , v_{d+1}$ of $S_{\sigma}$ are particularly simple, the denominator in Theorem \ref{theorem3} at a lattice point $y=m \in \mathbb{Z}^d$ simplifies as
\[ \left( z - \langle v_1 ,m \rangle \right) \cdots \left( z - \langle v_{d+1} ,m \rangle \right) = z \left( z - m_1 \sigma_1 a_1 \right) \cdots \left( z - m_d \sigma_d a_d \right) . \]
This means that the integrand in Theorem \ref{theorem3} can indeed have a high order pole at $z=0$, namely for lattice points $m \in \mathbb{Z}^d$ with many zero coordinates. We were able to find the sum of the residues at $z=0$ over all lattice points $m \in \mathbb{Z}^d$ and obtained the following.

\begin{definition}\label{definition2}
Let $a_1 , \dots , a_d >0$, and let $\zeta$ denote the Riemann zeta function. Let $p(t)=p_{(a_1 , \dots , a_d)}(t)=\sum_{k=0}^d c_k t^k$, where $c_d = \lambda (C) = \frac{2^d a_1 \cdots a_d}{d!}$, and
\[ c_k = \frac{2^d a_1 \cdots a_d}{(2 \pi i)^{d-k} k!} \sum_{\ell =1}^{d} \sum_{1 \le j_1 < \dots < j_{\ell} \le d} \sum_{\substack{i_1 + \cdots + i_{\ell}=d-k \\ i_1 , \dots , i_{\ell} \ge 2 \\ 2 \mid i_1 , \dots , i_{\ell}}} \frac{-2 \zeta (i_1)}{a_{j_1}^{i_1}} \cdots \frac{-2 \zeta (i_{\ell})}{a_{j_{\ell}}^{i_{\ell}}} \]
for $0 \le k \le d-1$.
\end{definition}

\noindent Let us also introduce a notation for the error terms, which come from the residues of simple poles at $z \neq 0$ of the integrand in Theorem \ref{theorem3}.

\begin{definition}\label{definition3}
Let $a_1 , \dots , a_d >0$, and let $N>0$ be an integer. Let
\[ E_N (t) = \sum_{j=1}^d \frac{i^d}{\pi^d N^d} \sum_{M \in [0,N-1]^d} \sum_{\substack{m \in [-M_1 , M_1] \times \cdots \times [-M_d,M_d] \\ m_j \neq 0}} \frac{e^{-2 \pi i m_j a_j t}}{m_j \prod_{k \neq j} \left( m_j \frac{a_j}{a_k} -m_k \right)} .  \]
\end{definition}

\noindent A combination of Theorem \ref{theorem2} and Theorem \ref{theorem3} thus yield the following.

\begin{prop}\label{proposition4}
Suppose that $\frac{1}{a_1} , \dots , \frac{1}{a_d} >0$ are algebraic and linearly independent over $\mathbb{Q}$. Let $C$ be as in \eqref{cross}. Then for any real $t>1$, any integer $N>1$ and any $\varepsilon >0$ we have
\[ \left| tC \cap \mathbb{Z}^d \right| = p(t) + E_N (t) + O \left( 1 + t^{d-1+ \varepsilon} \sqrt{\frac{\log N}{N}} \right) . \]
The implied constant depends only on $a_1 , \dots , a_d$ and $\varepsilon$, and is ineffective.
\end{prop}

\subsection{Statement of the main results}\label{section2.4}

The final step is to estimate the error terms $E_N (t)$, as in Definition \ref{definition3}. It is easy to see that the denominator in $E_N (t)$ is small, when the product
\[ \prod_{k \neq j} \left\| m_j \frac{a_j}{a_k} \right\| \]
is small. Thus we are interested in the following Diophantine quantity.

\begin{definition}\label{definition4}
For every integer $d \ge 1$ let $\gamma_d$ be the smallest real number $\gamma$ with the following property. If $\alpha_1 , \dots , \alpha_d$ are algebraic reals such that $1, \alpha_1 , \dots , \alpha_d$ are linearly independent over $\mathbb{Q}$, then
\[ \sum_{m=1}^M \frac{1}{\left\| m \alpha_1 \right\| \cdots \left\| m \alpha_d \right\|} = O \left( M^{\gamma + \varepsilon} \right) \]
for any $\varepsilon >0$ with an implied constant depending only on $\alpha_1 , \dots , \alpha_d$ and $\varepsilon$, as $M \to \infty$.
\end{definition}

It is easy to see that $1 \le \gamma_d \le 2$ for every $d$. Indeed, on the one hand, Dirichlet's theorem on Diophantine approximation states that there exist infinitely many positive integers $m$ such that
\[ \frac{1}{\left\| m \alpha_1 \right\| \cdots \left\| m \alpha_d \right\|} \ge \frac{1}{\left\| m \alpha_1 \right\|} \ge m , \]
which clearly shows $\gamma_d \ge 1$. On the other hand, applying Schmidt's theorem \eqref{Schmidt2} term by term we obtain $\gamma_d \le 2$.

A well-known argument based on the pigeonhole principle gives $\gamma_1 =1$. We were able to generalize that argument to higher dimensions to obtain the following result, which might be of interest in its own right.

\begin{thm}\label{theorem5}
For any $d \ge 1$ we have $\gamma_d \le 2-\frac{1}{d}$.
\end{thm}

Unfortunately we do not know if Theorem \ref{theorem5} is best possible for $d \ge 2$. In fact, we were not able to find any nontrivial lower bound for $\gamma_d$.

Our main result on the lattice point counting problem in the cross-polytope $C$ is the following. It is given in terms of the exponents $\gamma_d$ in the hope of future improvement on their values.

\begin{thm}\label{theorem6}
Suppose that $\frac{1}{a_1} , \dots , \frac{1}{a_d} >0$ are algebraic and linearly independent over $\mathbb{Q}$. Let $C, p(t)$ and $\gamma_d$ be as in \eqref{cross}, Definition \ref{definition2} and Definition \ref{definition4}.
\begin{itemize}
\item[(i)] For any $1 \le T_1 < T_2$ such that $T_2 - T_1 \ge 1$ we have
\[ \frac{1}{T_2 - T_1} \int_{T_1}^{T_2} \left( \left| tC \cap \mathbb{Z}^d \right| - p(t) \right) \, \mathrm{d}t = O(1) \]
with an ineffective implied constant depending only on $a_1 , \dots , a_d$.
\item[(ii)] For any real $t>1$ and $\varepsilon >0$ we have
\[ \left| tC \cap \mathbb{Z}^d \right| = p(t) + O \left( t^{\frac{\gamma_{d-1} -1}{\gamma_{d-1}} (d-1) + \varepsilon} \right) \]
with an ineffective implied constant depending only on $a_1 , \dots , a_d$ and $\varepsilon$.
\end{itemize}
\end{thm}

The lattice point counting problem in the simplex $S$, as in \eqref{simplex}, reduces to that in the cross-polytope $C$ using Proposition \ref{proposition1}. It is therefore natural to introduce the following polynomial.

\begin{definition}\label{definition5}
Let $a_1 , \dots , a_d >0$, and let $p_{(a_1 , \dots , a_d)}(t)$ be as in Definition \ref{definition2}. Let
\[ q(t) = q_{(a_1 , \dots , a_d)} (t) = \frac{1}{2^d} \sum_{I \subseteq [d]} p_{\left( a_i \,\, \middle| \,\, i \in I \right)} (t) . \]
\end{definition}

\noindent The main result on the lattice point counting problem in $S$ is thus the following.

\begin{thm}\label{theorem7}
Suppose that $\frac{1}{a_1} , \dots , \frac{1}{a_d} >0$ are algebraic and linearly independent over $\mathbb{Q}$. Let $S, q(t)$ and $\gamma_d$ be as in \eqref{simplex}, Definition \ref{definition5} and Definition \ref{definition4}.
\begin{itemize}
\item[(i)] For any $1 \le T_1 < T_2$ such that $T_2 - T_1 \ge 1$ we have
\[ \frac{1}{T_2 - T_1} \int_{T_1}^{T_2} \left( \left| tS \cap \mathbb{Z}^d \right| - q(t) \right) \, \mathrm{d}t = O(1) \]
with an ineffective implied constant depending only on $a_1 , \dots , a_d$.
\item[(ii)] For any real $t>1$ and $\varepsilon >0$ we have
\[ \left| tS \cap \mathbb{Z}^d \right| = q(t) + O \left( t^{\frac{\gamma_{d-1} -1}{\gamma_{d-1}} (d-1) + \varepsilon} \right) \]
with an ineffective implied constant depending only on $a_1 , \dots , a_d$ and $\varepsilon$.
\end{itemize}
\end{thm}

Theorems \ref{theorem6} (ii) and \ref{theorem7} (ii) were stated in terms of the unknown quantity $\gamma_d$. The estimate in Theorem \ref{theorem5} gives the following bounds.

\begin{thm}\label{theorem8}
Suppose that $\frac{1}{a_1} , \dots , \frac{1}{a_d} >0$ are algebraic and linearly independent over $\mathbb{Q}$. Let $C$, $S$, $p(t)$ and $q(t)$ be as in \eqref{cross}, \eqref{simplex}, Definition \ref{definition2} and Definition \ref{definition5}. For any real $t>1$ and $\varepsilon >0$ we have
\[ \begin{split} \left| tC \cap \mathbb{Z}^d \right| &= p(t) + O \left( t^{\frac{(d-1)(d-2)}{2d-3} + \varepsilon} \right) , \\ \left| tS \cap \mathbb{Z}^d \right| &= q(t) + O \left( t^{\frac{(d-1)(d-2)}{2d-3} + \varepsilon} \right) \end{split} \]
with ineffective implied constants depending only on $a_1 , \dots , a_d$ and $\varepsilon$.
\end{thm}
\begin{flushright}
$\square$
\end{flushright}

\subsection{Conclusions}\label{section2.5}

Let us now list some corollaries and remarks on the main results.
\begin{enumerate}
\item Theorems \ref{theorem6} (i), \ref{theorem7} (i) clearly show that $p(t)$ and $q(t)$ are indeed the main terms of $\left| tC \cap \mathbb{Z}^d \right|$ and $\left| tS \cap \mathbb{Z}^d \right|$, respectively. This means that our intuition about the residues of the high order poles in Theorem \ref{theorem3} being the main contribution in the Poisson summation formula was correct.

Several examples of compact sets $B \subset \mathbb{R}^d$ are known for which the number of lattice points $\left| tB \cap \mathbb{Z}^d \right|$, as a function of the real variable $t>1$ can be approximated by a function other than the Lebesgue measure $\lambda (B) t^d$. Let us only mention the example of the torus
\[ B = \left\{ (x,y,z) \in \mathbb{R}^3 \,\, \middle| \,\, \left( \sqrt{x^2 + y^2} - a \right)^2 + z^2 \le b^2 \right\} , \]
where $0<b<a$ are constants. Nowak \cite{Nowak} proves
\[ \left| tB \cap \mathbb{Z}^3 \right| = \lambda (B) t^3 + F_{a,b} (t) t^{\frac{3}{2}} + O \left( t^{\frac{11}{8} + \varepsilon} \right) \]
for any $\varepsilon >0$, where $F_{a,b}$ is a bounded function defined by the absolutely convergent trigonometric series
\[ F_{a,b} (t) = 4 a \sqrt{b} \sum_{n=1}^{\infty} n^{- \frac{3}{2}} \sin \left( 2 \pi n b t - \frac{\pi}{4} \right) . \]
Here the second order term $F_{a,b}(t) t^{\frac{3}{2}}$ is related to the points on the boundary $\partial B$ with Gaussian curvature zero.

\item Theorem \ref{theorem8} in dimension $d=2$ gives the error bound $O \left( t^{\varepsilon} \right)$ of Skriganov \cite{Skrig2}. Any improvement on Theorem \ref{theorem5} would result in better error bounds in higher dimensions. E.g. if $\gamma_{d-1}=1$, then the error is $O \left( t^{\varepsilon} \right)$ in dimension $d$.

\item Even though we allowed the dilation factor $t$ to be a real number, the main terms $p(t)$ and $q(t)$ were polynomials. In contrast, for a rational polytope $P \subset \mathbb{R}^d$, $|tP \cap \mathbb{Z}^d|$ is a quasi-polynomial, but not a polynomial as a function of the \textit{real} variable $t$. It is thus more natural to compare our polynomials $p(t)$ and $q(t)$ to Ehrhart polynomials, defined via \textit{integral} dilations of a lattice polytope. Despite the fact that their natural domains are different, $p(t)$ and $q(t)$ seem to show a certain similarity to Ehrhart polynomials. Without providing a deeper understanding, let us mention a few of these similarities.

Definition \ref{definition2} of $p(t) = \sum_{k=0}^d c_k t^k$ gives that for any $k \not\equiv d \pmod{2}$ we have $c_k =0$. Indeed, for such $k$ the number $d-k$ cannot be written as a sum of positive even integers, resulting in an empty sum defining $c_k$. In other words, the polynomial $p(t)$ satisfies the functional equation $p(-t)=(-1)^d p(t)$. Note that for any lattice polytope $P$ there exists a polynomial $f(t)$ such that
\[ f(t) = \left| tP \cap \mathbb{Z}^d \right| - \frac{1}{2} \left| t \left( \partial P \right) \cap \mathbb{Z}^d \right| \]
for every positive integer $t$, and that this polynomial also satisfies the functional equation $f(-t)=(-1)^d f(t)$. This is a form of the famous \textbf{Ehrhart--Macdonald reciprocity} \cite{Macd}. This shows a clear connection between $p(t)$ and Ehrhart polynomials, even though $C$ is not a lattice polytope.

\item In Definition \ref{definition2} of the coefficients $c_k$ of $p(t)$ we have
\[ \zeta (i_1) \cdots \zeta (i_{\ell}) \in \pi^{i_1 + \cdots + i_{\ell}} \mathbb{Q} = \pi^{d-k} \mathbb{Q} , \]
therefore $c_k$ is a rational function of $a_1, \dots , a_d$ with rational coefficients. The first two nontrivial coefficients are
\[ \begin{split} c_{d-2} &= \frac{2^{d-2} a_1 \cdots a_d}{3 (d-2)!} \sum_{1 \le i \le d} \frac{1}{a_i^2} , \\ c_{d-4} &= \frac{2^{d-4} a_1 \cdots a_d}{9 (d-4)!} \left( \sum_{1 \le i < j \le d} \frac{1}{a_i^2 a_j^2} - \frac{1}{5} \sum_{1 \le i \le d} \frac{1}{a_i^4} \right) . \end{split} \]

In particular, $c_{d-2}>0$. Under the assumptions of Theorem \ref{theorem6} the coordinates of every normal vector of $C$ are algebraic and linearly independent over $\mathbb{Q}$, yet the lattice discrepancy satisfies
\[ \left| tC \cap \mathbb{Z}^d \right| - \lambda (C)t^d \sim c_{d-2} t^{d-2} . \]
This shows that Randol's theorem \cite{Rand} mentioned in the Introduction is best possible even under stronger conditions.

\item Definitions \ref{definition2}, \ref{definition4} show that the coefficients of $q(t)$ are also rational functions of $a_1 , \dots , a_d$ with rational coefficients. Writing $q(t) = \sum_{k=0}^d e_k t^k$ we clearly have $e_d = \lambda (S) = \frac{a_1 \cdots a_d}{d!}$. The next few coefficients are
\[ \begin{split} e_{d-1} &= \frac{a_1 \cdots a_d}{2 (d-1)!} \sum_{1 \le i \le d} \frac{1}{a_i} , \\ e_{d-2} &= \frac{a_1 \cdots a_d}{4 (d-2)!} \left( \frac{1}{3} \sum_{1 \le i \le d} \frac{1}{a_i^2} + \sum_{1 \le i<j \le d} \frac{1}{a_i a_j} \right) , \\ e_{d-3} &= \frac{a_1 \cdots a_d}{8 (d-3)!} \left( \frac{1}{3} \sum_{1 \le i<j \le d} \left( \frac{1}{a_i a_j^2} + \frac{1}{a_i^2 a_j} \right) + \sum_{1 \le i<j<k \le d} \frac{1}{a_i a_j a_k} \right) . \end{split} \]
Note that $e_{d-1}$ is one half of the total surface area of the $d-1$ dimensional faces of $S$ with a rational equation. This is perfect analogy with Ehrhart polynomials, if we use the natural convention that the ``sublattice'' of $\mathbb{Z}^d$ on the affine hyperplane with normal vector $\left( \frac{1}{a_1} , \dots , \frac{1}{a_d} \right)$ (in fact the empty set or a singleton) has infinite covolume, making the normalized surface area of the face zero.

In the case when $a_1 , \dots , a_d$ are positive \textit{integers}, the simplex $S$ has an actual Ehrhart polynomial. This Ehrhart polynomial has been computed using methods as diverse as the theory of toric varieties \cite{Pomm}, Fourier analysis \cite{Diaz} and complex analysis \cite{MBeck}. If $a_1 , \dots , a_d$ are pairwise coprime integers, the coefficient of $t^{d-2}$ in this Ehrhart polynomial is
\begin{multline*}\frac{a_1 \cdots a_d}{4 (d-2)!} \left( \frac{1}{3} \sum_{1 \le i \le d} \frac{1}{a_i^2} + \sum_{1 \le i<j \le d} \frac{1}{a_i a_j} \right) \\ + \frac{1}{(d-2)!} \left( \frac{d}{4} + \frac{1}{12 a_1 \cdots a_d} - \sum_{1 \le i \le d} s \left( \frac{a_1 \cdots a_d}{a_i} , a_i \right) \right) , \end{multline*}
where $s$ is the Dedekind sum defined as
\[ s(a,b) = \sum_{k=1}^{b-1} \left( \frac{k}{b} - \frac{1}{2} \right) \left( \left\{ \frac{ak}{b} \right\} - \frac{1}{2} \right) \]
for coprime integers $a,b$.
\end{enumerate}

\section{Proofs}\label{section3}

In this Section we give the proofs of the results in the same order in which they were stated.

\vspace{5mm}

\noindent\textbf{Proof of Theorem \ref{theorem2}:} We start with the following lemma, which will help estimate the number of lattice points close to the boundary of $tP$.

\begin{lemma}\label{lemma9} Let $2 \le k \le d$, and suppose that the coordinates of $n= (n_1 , \dots , n_d )$ are algebraic and span a vector space of dimension $k$ over $\mathbb{Q}$. Let $B \subset \mathbb{R}^d$ be a ball of radius $R>1$, and consider two parallel affine hyperplanes orthogonal to $n$ at distance $a>0$ from each other. Then the number of lattice points in $B$ which fall between the two affine hyperplanes is $O \left( R^{d-k} + a R^{d-1+ \varepsilon} \right)$ for any $\varepsilon >0$. The implied constant depends only on $n$ and $\varepsilon$, and is ineffective.
\end{lemma}

\noindent\textbf{Proof of Lemma \ref{lemma9}:} We may assume $n_d=1$. The region we are interested in is
\[ A = \left\{ x \in B \,\, \middle| \,\, b \le \left\langle \frac{n}{|n|} , x \right\rangle \le b+a \right\} \]
for some $b \in \mathbb{R}$.

Let $\alpha_1 , \dots , \alpha_{k-1} , \alpha_k$ be a basis in the vector space spanned by $n_1 , \dots , n_d$ over $\mathbb{Q}$, such that $\alpha_k=1$. Schmidt's theorem \eqref{Schmidt1} states that
\begin{equation}\label{eq3}
\left\| m_1 \alpha_1 + \cdots + m_{k-1} \alpha_{k-1} \right\| \ge \frac{K}{|m|^{k-1+ \varepsilon}}
\end{equation}
for any $m \in \mathbb{Z}^{k-1} \backslash \{ 0 \}$, with some constant $K>0$ depending only on $\alpha_1 , \dots , \alpha_{k-1}$ and $\varepsilon$. Since $\alpha_1 , \dots , \alpha_k$ is a basis, we have
\[ n_i = \sum_{j=1}^k \frac{\lambda_{i,j}}{Q} \alpha_j \]
for some $\lambda_{i,j} \in \mathbb{Z}$ and $Q \in \mathbb{N}$.

Let $c,c' \in A \cap \mathbb{Z}^d$ be such that $\langle c-c' , n \rangle \neq 0$. Then
\[ \left| \left\langle c-c' , \frac{n}{|n|} \right\rangle \right| = \frac{1}{Q |n|} \left| \sum_{j=1}^k \sum_{i=1}^d (c_i-c_i') \lambda_{i,j} \alpha_j \right| \ge \frac{1}{Q |n|} \left\| \sum_{j=1}^{k-1} \sum_{i=1}^d (c_i-c_i') \lambda_{i,j} \alpha_j \right\| , \]
since the $j=k$ term is an integer. Let $m_j = \sum_{i=1}^d (c_i-c_i') \lambda_{i,j} \in \mathbb{Z}$ for $1 \le j \le k-1$. If $m \in \mathbb{Z}^{k-1} \backslash \{ 0 \}$, then \eqref{eq3} implies
\[ \left| \left\langle c-c' , \frac{n}{|n|} \right\rangle \right| \ge \frac{K'}{|m|^{k-1+\varepsilon}} \]
for some $K'>0$. Clearly $|m| = O \left( |c-c'| \right)$. Since $c,c'$ lie in a ball of radius $R$ we obtain
\begin{equation}\label{eq4}
\left| \left\langle c-c' , \frac{n}{|n|} \right\rangle \right| \ge \frac{K''}{R^{k-1+\varepsilon}}
\end{equation}
for some $K''>0$. \eqref{eq4} is clearly true in the case $m=0$ as well.

The geometric meaning of \eqref{eq4} is the following. Let us draw an affine hyperplane with normal vector $n$ through every lattice point $c \in A \cap \mathbb{Z}^d$. Then the distance of any two of these hyperplanes is at least $\frac{K''}{R^{k-1+ \varepsilon}}$. Hence the number of such hyperplanes is $O \left( \lceil a R^{k-1+\varepsilon} \rceil \right)$. Every such hyperplane contains a sublattice of $\mathbb{Z}^d$ of dimension $d-k$. Therefore the number of lattice points on a given hyperplane inside $B$ is $O \left( R^{d-k} \right)$. The total number of lattice points in $A$ is thus
\[ O \left( \lceil a R^{k-1+ \varepsilon } \rceil R^{d-k} \right) = O \left( R^{d-k} + a R^{d-1 + \varepsilon} \right) . \]
\begin{flushright}
$\square$
\end{flushright}

The Fej\'er kernel corresponding to the Ces\`aro means in Definition \ref{definition1} is the function $F_N : \mathbb{R}^d \to \mathbb{R}$ defined as
\[ F_N (x) = \frac{1}{N^d} \sum_{M \in [0,N-1]^d} \sum_{m \in [-M_1 , M_1] \times \cdots \times [-M_d,M_d]} e^{2 \pi i \langle m , x \rangle} . \]
For the basic properties of $F_N$ see e.g.\ Section 3.1.3.\ in \cite{Graf}. Introducing the function $f: \left[ - \frac{1}{2} , \frac{1}{2} \right]^d \to \mathbb{R}$ defined as
\[ f(x) = \sum_{m \in \mathbb{Z}^d} \chi_{tP} (m+x) , \]
we have that
\begin{equation}\label{integral}
\mathrm{Ces} (tP , N) - L = \int_{\left[ - \frac{1}{2} , \frac{1}{2} \right]^d} \left( f(x) -L \right) F_N (x) \, \mathrm{d}x
\end{equation}
for any $L \in \mathbb{R}$. In the $d=1$ case it is well known that $F_N \ge 0$ and that for any $0<h<\frac{1}{2}$ we have
\[ \int_{\left[ - \frac{1}{2} , \frac{1}{2} \right] \backslash [-h,h]} F_N (x) \, \mathrm{d}x = O \left( \frac{\log N}{h N} \right) , \]
the latter being an easy exercise using summation by parts. Since the $d$ dimensional Fej\'er kernel factors into one dimensional ones as $F_N (x_1 , \dots , x_d) = F_N (x_1) \cdots F_N (x_d)$, we obtain that $F_N \ge 0$ holds in any dimension. Recalling that the total integral of $F_N$ over $\left[-\frac{1}{2}, \frac{1}{2} \right]^d$ is 1, Fubini's theorem implies that
\begin{equation}\label{Fejer}
\int_{\left[ - \frac{1}{2} , \frac{1}{2} \right]^d \backslash [-h,h]^d} F_N (x) \, \mathrm{d}x = O \left( \frac{\log N}{h N} \right)
\end{equation}
holds for any $0<h<\frac{1}{2}$ in any dimension as well, with an implied constant depending only on $d$.

Let $0<h<\frac{1}{2}$ be arbitrary, and use \eqref{integral} with $L=\left| tP \cap \mathbb{Z}^d \right|$ to get
\begin{equation}\label{ineq}\begin{split} \left| \mathrm{Ces} (tP ,N) - \left| tP \cap \mathbb{Z}^d \right| \right| \le &\int_{[-h,h]^d} \left| f(x) - \left| tP \cap \mathbb{Z}^d \right| \right| F_N (x) \, \mathrm{d}x \\ &+ \int_{\left[ - \frac{1}{2} , \frac{1}{2} \right]^d \backslash [-h,h]^d} \left| f(x) - \left| tP \cap \mathbb{Z}^d \right| \right| F_N (x) \, \mathrm{d}x . \end{split}\end{equation}

To estimate the first integral in \eqref{ineq} note that for any $x \in [-h,h]^d$ we have
\[ \begin{split} \left| f(x) - \left| tP \cap \mathbb{Z}^d \right| \right| &\le \sum_{m \in \mathbb{Z}^d} \left| \chi_{tP} (m+x) - \chi_{tP} (m) \right| \\ &\le \left| \left\{ m \in \mathbb{Z}^d \,\, \middle| \,\, \mathrm{dist} \left( m , \partial (tP) \right) \le \sqrt{d} h \right\} \right| , \end{split}\]
where $\mathrm{dist} (y,A)$ denotes the distance of a point $y \in \mathbb{R}^d$ from a set $A \subseteq \mathbb{R}^d$. The set
\[ \left\{ y \in \mathbb{R}^d \,\, \middle| \,\, \mathrm{dist} \left( y , \partial (tP) \right) \le \sqrt{d} h \right\} \]
can be covered by regions as in Lemma \ref{lemma9} with $R= O(t)$ and $a=O(h)$. Moreover, the number of such regions required is the number of $d-1$ dimensional faces of $P$. Thus
\[ \left| f(x) - \left| tP \cap \mathbb{Z}^d \right| \right| = O \left( t^{d-k} + h t^{d-1+ \varepsilon} \right) \]
for any $x \in [-h,h]^d$, and hence
\begin{equation}\label{int1}
\int_{[-h,h]^d} \left| f(x) - \left| tP \cap \mathbb{Z}^d \right| \right| F_N (x) \, \mathrm{d}x = O \left( t^{d-k} + h t^{d-1+ \varepsilon} \right) .
\end{equation}

It is not difficult to see that the error term in \eqref{eq2} is invariant under translations of the polytope. In other words, we have the slightly more general estimate
\[ \left| (tP -x) \cap \mathbb{Z}^d \right| = \lambda (P) t^d + O(t^{d-1}) \]
for any $x \in \mathbb{R}^d$, with an implied constant depending only on $P$ but not on $x$. In the second integral of \eqref{ineq} we thus have
\[ \left| f(x) - \left| tP \cap \mathbb{Z}^d \right| \right| = \left| \left| (tP-x) \cap \mathbb{Z}^d \right| - \left| tP \cap \mathbb{Z}^d \right| \right| = O \left( t^{d-1} \right) \]
with an implied constant independent of $x$. Therefore \eqref{Fejer} implies
\begin{equation}\label{int2}
\int_{\left[ - \frac{1}{2} , \frac{1}{2} \right]^d \backslash [-h,h]^d} \left| f(x) - \left| tP \cap \mathbb{Z}^d \right| \right| F_N (x) \, \mathrm{d}x = O \left( t^{d-1} \frac{\log N}{hN} \right) .
\end{equation}

Using \eqref{ineq}, \eqref{int1} and \eqref{int2} we obtain
\[ \mathrm{Ces} (tP ,N) - \left| tP \cap \mathbb{Z}^d \right| = O \left( t^{d-k} + h t^{d-1+ \varepsilon} + t^{d-1} \frac{\log N}{hN} \right) \]
for any $0<h<\frac{1}{2}$. Choosing $h= \sqrt{\frac{\log N}{N}}$ to minimize the error finishes the proof of Theorem \ref{theorem2}.
\begin{flushright}
$\square$
\end{flushright}

\noindent\textbf{Proof of Theorem \ref{theorem3}:} Consider the simplex
\[ S_0 = \left\{ x \in \mathbb{R}^d \,\, \middle| \,\, x_1 , \dots , x_d \ge 0 , \,\, x_1 + \cdots + x_d \le 1 \right\} , \]
let $t>0$ be real, and let $y \in \mathbb{R}^d$ be such that $y_j \neq 0$ and $y_j \neq y_k$ for any $j \neq k$. We shall prove that
\begin{equation}\label{induction}
\hat{\chi}_{tS_0} (y) = \frac{(-1)^{d+1}}{(2 \pi i)^d} \sum_{j=1}^d \frac{1-e^{-2 \pi i y_j t}}{y_j \prod_{k \neq j} (y_j-y_k)}
\end{equation}
by induction on $d$. The $d=1$ case is trivial, using the convention that an empty product is 1. Suppose the claim holds in dimension $d-1$, fix $x_d \in [0,t]$ and consider the cross section
\begin{multline*} \left\{ (x_1 , \dots , x_{d-1}) \in \mathbb{R}^{d-1} \,\, \middle| \,\, (x_1 , \dots , x_d) \in tS_0 \right\} \\= \left\{ (x_1 , \dots , x_{d-1}) \in \mathbb{R}^{d-1} \,\, \middle| \,\, x_1, \dots , x_{d-1} \ge 0 \,\, x_1 + \cdots + x_{d-1} \le t-x_d \right\} . \end{multline*}
The inductive hypothesis with $t-x_d$ instead of $t$, and Fubini's theorem thus imply that
\[ \begin{split} \hat{\chi}_{tS_0}(y) &= \int_0^t \frac{(-1)^d}{(2 \pi i)^{d-1}} \sum_{j=1}^{d-1} \frac{1-e^{-2 \pi i y_j (t-x_d)}}{y_j \prod_{k \neq j,d} (y_j - y_k)} e^{-2 \pi i y_d x_d} \, \mathrm{d}x_d \\ &= \frac{(-1)^{d+1}}{(2 \pi i)^d} \sum_{j=1}^{d-1} \frac{1-e^{-2 \pi i y_j t}}{y_j \prod_{k \neq j} (y_j - y_k)} \\ & \hspace{5mm} + \frac{(-1)^{d+1}}{(2 \pi i)^d} \left( \sum_{j=1}^{d-1} \frac{-1}{y_d \prod_{k \neq j} (y_j - y_k)} \right) (1-e^{-2 \pi i y_d t}) . \end{split} \]
To finish the proof of \eqref{induction} we need to show
\begin{equation}\label{finish}
\sum_{j=1}^{d-1} \frac{-1}{y_d \prod_{k \neq j} (y_j - y_k)} = \frac{1}{y_d \prod_{k \neq d} (y_d - y_k)} .
\end{equation}
To this end, consider the partial fraction decomposition
\begin{equation}\label{partial}
\frac{1}{\prod_{k=1}^{d-1} (x-y_k)} = \sum_{j=1}^{d-1} \frac{A_j}{x-y_j} ,
\end{equation}
where the constant $A_j$ is
\[ A_j = \frac{1}{\prod_{k \neq j,d} (y_j-y_k)} . \]
Substituting $x=y_d$ in \eqref{partial} we obtain \eqref{finish}, which in turn finishes the proof of \eqref{induction}.

The main idea is to identify the formula found in \eqref{induction} as the sum of residues of a meromorphic function. For any $y \in \mathbb{R}^d$ such that $y_j \neq 0$ and $y_j \neq y_k$ for any $j \neq k$ we have
\[ \frac{(-1)^{d+1}}{(2 \pi i)^d} \sum_{j=1}^d \frac{1-e^{-2 \pi i y_j t}}{y_j \prod_{k \neq j} (y_j-y_k)} = \frac{(-1)^{d+1}}{(2 \pi i)^{d+1}} \int_{|z|=R} \frac{1-e^{-2 \pi i z t}}{z (z-y_1) \cdots (z-y_d)} \, \mathrm{d}z \]
for any $R> \max_j |y_j|$. Indeed, the meromorphic integrand has $d+1$ distinct isolated singularities. The singularity at $z=0$ is removable, while the singularity at $z=y_j$ is a simple pole the residue of which is exactly the $j$th term of the sum.

We now claim that
\begin{equation}\label{eq5}
\hat{\chi}_{tS_0} (y) = \frac{(-1)^{d+1}}{(2 \pi i)^{d+1}} \int_{|z|=R} \frac{1-e^{-2 \pi i z t}}{z (z-y_1) \cdots (z-y_d)} \, \mathrm{d}z
\end{equation}
holds for \textit{any} $y \in \mathbb{R}^d$, as long as $R> \max_j |y_j|$. Fix an arbitrary constant $r>0$. It is enough to show \eqref{eq5} in the ball $|y| \le r$. From the definition of the Fourier transform and Lebesgue's dominated convergence theorem we get that the left hand side of \eqref{eq5} is a continuous function of $y$. It is easy to see that the right hand side of \eqref{eq5} is also a continuous function of $y$ on the ball $|y| \le r$, by choosing $R>r$. Since these continuous functions are equal on a dense subset of the ball $|y| \le r$, they are equal everywhere.

Note that
\[ \int_{|z|=R} \frac{1}{z (z-y_1) \cdots (z-y_d)} \, \mathrm{d}z =0 \]
for $R> \max_j |y_j|$. Indeed, the residue theorem implies that the value of the integral does not depend on $R$. On the other hand, the trivial estimate gives that the integral is $O \left( R^{-d} \right)$, as $R \to \infty$. Therefore
\begin{equation}\label{eq6}
\hat{\chi}_{tS_0} (y) = \frac{(-1)^d}{(2 \pi i)^{d+1}} \int_{|z|=R} \frac{e^{-2 \pi i z t}}{z (z-y_1) \cdots (z-y_d)} \, \mathrm{d}z
\end{equation}
for any $R > \max_j |y_j|$.

Now let $S \subset \mathbb{R}^d$ be an arbitrary simplex with vertices $v_1, \dots , v_{d+1}$. Let $M$ be the $n \times n$ matrix the columns of which are the vectors $v_1 - v_{d+1} , \dots , v_d - v_{d+1}$, and let $g(x) = Mx + t v_{d+1}$. Then $g(tS_0)=tS$, thus using $g(x)$ as an integral transformation we get
\begin{equation}\label{eq7}
\hat{\chi}_{tS} (y) = \int_{tS_0} e^{- 2 \pi i \langle Mx+ t v_{d+1},y \rangle} \left| \det M \right| \, \mathrm{d}x .
\end{equation}
Since $\lambda (S_0) = \frac{1}{d!}$, substituting $t=1$ and $y=0$ in \eqref{eq7} we obtain $\left| \det M \right| = d! \lambda (S)$. Therefore \eqref{eq7} yields
\[ \hat{\chi}_{tS} (y) = d! \lambda (S) e^{- 2 \pi i \langle v_{d+1} ,y \rangle t} \hat{\chi}_{tS_0} \left( M^T y \right) , \]
where $M^T$ denotes the transpose of $M$. The coordinates of the vector $M^T y$ are
\[ \langle v_1 - v_{d+1} , y \rangle , \dots , \langle v_d - v_{d+1} , y \rangle , \]
hence \eqref{eq6} gives
\[ \hat{\chi}_{tS} (y) = \frac{(-1)^d d!}{(2 \pi i)^{d+1}} \lambda (S) \int_{|z|=R} \frac{e^{- 2 \pi i (z + \langle v_{d+1} ,y \rangle ) t}}{z \left( z - \langle v_1 - v_{d+1} , y \rangle \right) \cdots \left( z - \langle v_d - v_{d+1} , y \rangle \right) } \, \mathrm{d}z , \]
where $R > \max_j \left| \langle v_j - v_{d+1} , y \rangle \right|$. Finally, let us apply the simple integral transformation $f(z)=z- \langle v_{d+1} ,y \rangle$, to get
\[ \hat{\chi}_{tS} (y) = \frac{(-1)^d d!}{(2 \pi i)^{d+1}} \lambda (S) \int_{\gamma} \frac{e^{- 2 \pi i z t}}{\left( z - \langle v_1 , y \rangle \right) \cdots \left( z - \langle v_{d+1} , y \rangle \right) } \, \mathrm{d}z , \]
where $\gamma$ is a circle centered at $\langle v_{d+1} , y \rangle$ which contains every singularity of the integrand inside. The residue theorem implies that we can replace $\gamma$ by a circle centered at the origin of radius $R> \max_j \left| \langle v_j, y \rangle \right|$.
\begin{flushright}
$\square$
\end{flushright}

\noindent\textbf{Proof of Proposition \ref{proposition4}:} Theorem \ref{theorem2} implies that
\begin{equation}\label{final}
\left| tC \cap \mathbb{Z}^d \right| = \mathrm{Ces} (tC , N) + O \left( 1+ t^{d-1+\varepsilon } \sqrt{\frac{\log N}{N}} \right) ,
\end{equation}
where $\mathrm{Ces} (tC , N)$ is as in Definition \ref{definition1}. The simplices $S_{\sigma}$, as in \eqref{triangulation}, $\sigma \in \left\{ 1, -1 \right\}^d$,  triangulate $C$, therefore
\[ \hat{\chi}_{tC} = \sum_{\sigma \in \left\{ 1, -1 \right\}^d} \hat{\chi}_{tS_{\sigma}} . \]
It is easy to see that
\[ \mathrm{Ces} (tC, N) = \sum_{\sigma \in \left\{ 1, -1 \right\}^d} \mathrm{Ces} (tS_{\sigma}, N) = 2^d \, \mathrm{Ces} (tS, N) , \]
where $S$ is as in \eqref{simplex}. Applying Theorem \ref{theorem3} to $S$ with a fixed $R> N \max_{j} a_j$, and substituting $\lambda (S) = \frac{a_1 \cdots a_d}{d!}$ we obtain
\begin{equation}\label{eq8}
\mathrm{Ces} (tC, N) = \frac{1}{N^d} \sum_{M \in [0,N-1]^d} A_M
\end{equation}
with
\begin{equation}\label{AM}
A_M = \sum_{m \in [-M_1 , M_1] \times \cdots \times [-M_d, M_d]} \frac{(-1)^d 2^d a_1 \cdots a_d}{(2 \pi i)^{d+1}} \int_{|z|=R} \frac{e^{-2 \pi i z t}}{z (z - m_1 a_1) \cdots (z- m_d a_d)} \, \mathrm{d}z .
\end{equation}

We now wish to apply the residue theorem to the complex line integral in \eqref{AM}. Note that the pole at $m_j a_j$ for $m_j \neq 0$ is simple. To separate the residue of the pole at $z=0$ from that of other poles, let us introduce
\[ B_M = \sum_{m \in [-M_1 , M_1] \times \cdots \times [-M_d, M_d]} \frac{(-1)^d 2^d a_1 \cdots a_d}{(2 \pi i)^d} \mathrm{Res}_0 \frac{e^{-2 \pi i z t}}{z (z - m_1 a_1) \cdots (z- m_d a_d)} . \]
Recalling Definition \ref{definition3}, \eqref{eq8} hence simplifies as
\begin{equation}\label{Ces}
\mathrm{Ces}(tC,N) = \frac{1}{N^d} \sum_{M \in [0,N-1]^d} B_M + E_N (t) .
\end{equation}

It is easy to see that if $m=0$, then the residue in question is
\[ \frac{(-1)^d 2^d a_1 \cdots a_d}{(2 \pi i)^d} \mathrm{Res}_0 \frac{e^{-2 \pi i z t}}{z^{d+1}} = \lambda (C) t^d . \]
Let us now fix a lattice point $m \in \mathbb{Z}^d \backslash \{ 0 \}$. Suppose $m$ has exactly $\ell$ nonzero coordinates, $m_{j_1} , \dots , m_{j_{\ell}} \neq 0$, for some $1 \le \ell \le d$ and $1 \le j_1 < \cdots < j_{\ell} \le d $. Using well-known Taylor series expansions we obtain that
\[ \mathrm{Res}_0 \frac{e^{-2 \pi i z t}}{z (z - m_1 a_1) \cdots (z- m_d a_d)} = \mathrm{Res}_0 \frac{1}{z^{d+1}} e^{-2 \pi i z t} \frac{z}{z-m_{j_1} a_{j_1}} \cdots \frac{z}{z-m_{j_{\ell}} a_{j_{\ell}}} \]
equals the coefficient of $z^d$ in the power series
\[ \left( \sum_{k=0}^{\infty} \frac{(-2 \pi i t)^k}{k!} z^k \right) \left( \sum_{i_1=1}^{\infty} \frac{-1}{(m_{j_1} a_{j_1})^{i_1}} z^{i_1} \right) \cdots \left( \sum_{i_{\ell}=1}^{\infty} \frac{-1}{(m_{j_{\ell}} a_{j_{\ell}})^{i_{\ell}}} z^{i_{\ell}} \right) . \]
Hence for such an $m$ we have
\begin{multline}\label{residues}
\frac{(-1)^d 2^d a_1 \cdots a_d}{(2 \pi i)^d} \mathrm{Res}_0 \frac{e^{-2 \pi i z t}}{z (z - m_1 a_1) \cdots (z- m_d a_d)} \\ = \sum_{k=0}^{d-1} \frac{2^d a_1 \cdots a_d}{(-2 \pi i)^{d-k}k!} t^k \sum_{\substack{i_1 + \cdots + i_{\ell} = d-k \\ i_1 , \dots , i_{\ell} \ge 1 }} \frac{-1}{(m_{j_1} a_{j_1})^{i_1}} \cdots \frac{-1}{(m_{j_{\ell}} a_{j_{\ell}})^{i_{\ell}}} .
\end{multline}
The sum of \eqref{residues} over $m_{j_1} \in [-M_{j_1} , M_{j_1}] \backslash \{ 0 \}, \dots , m_{j_{\ell}} \in [-M_{j_{\ell}} , M_{j_{\ell}}] \backslash \{ 0 \}$ is clearly
\[ \sum_{k=0}^{d-2} \frac{2^d a_1 \cdots a_d}{(2 \pi i)^{d-k}k!} t^k \sum_{\substack{i_1 + \cdots + i_{\ell} = d-k \\ i_1 , \dots , i_{\ell} \ge 2 \\ 2 \mid i_1 , \dots , i_{\ell} }} \frac{-2 \zeta (i_1)}{a_{j_1}^{i_1}} \cdots \frac{-2 \zeta (i_{\ell})}{ a_{j_{\ell}}^{i_{\ell}}} + O \left( \frac{t^{d-2}}{M_{j_1} +1} + \cdots + \frac{t^{d-2}}{M_{j_{\ell}}+1} \right) . \]
Recalling Definition \ref{definition2} we thus obtain
\begin{equation}\label{BM}
\begin{split}B_M &= p(t) + O \left( \frac{t^{d-2}}{M_1 +1} + \cdots + \frac{t^{d-2}}{M_d +1} \right) , \\ \frac{1}{N^d} \sum_{M \in [0,N-1]^d} B_M &= p(t) + O \left( t^{d-2} \frac{\log N}{N} \right) .\end{split}
\end{equation}
Combining \eqref{final}, \eqref{Ces} and \eqref{BM} concludes the proof.
\begin{flushright}
$\square$
\end{flushright}

\noindent\textbf{Proof of Theorem \ref{theorem5}:} Given irrational numbers $\alpha_1 , \dots , \alpha_d$, let
\begin{equation}\label{lower}
L_M = \min_{1 \le m \le M} \left\| m \alpha_1 \right\| \cdots \left\| m \alpha_d \right\|
\end{equation}
for any positive integer $M$. Clearly $0<L_M< \frac{1}{2^d}$. For any real number $h>1$ consider the set
\[ A_h = \left\{ 1 \le m \le M \,\, \middle| \,\, \left\| m \alpha_1 \right\| \cdots \left\| m \alpha_d \right\| < h L_M \right\} . \]
We wish to find an upper bound to the cardinality of $A_h$.

For any real number $0<c< \frac{1}{2^d}$ consider the set
\[ U_c = \left\{ x \in \left[ - \frac{1}{2} , \frac{1}{2} \right)^d \,\, \middle| \,\, \left| x_1 \cdots x_d \right| < c \right\} . \]
We shall prove by induction on $d$ that $\lambda (U_c) = O \left( c \log^{d-1} \frac{1}{c} \right) $ with an implied constant depending only on $d$. The case $d=1$ is trivial. Suppose the claim holds in dimension $d-1$. Fix an arbitrary $x_d \in \left[ - \frac{1}{2} , \frac{1}{2} \right) \backslash \{ 0 \}$, and consider the cross section
\begin{multline*}\left\{ (x_1 , \dots , x_{d-1}) \in \left[ - \frac{1}{2} , \frac{1}{2} \right)^{d-1} \,\, \middle| \,\, (x_1 , \dots , x_d) \in U_c \right\} \\ = \left\{ (x_1 , \dots , x_{d-1}) \in \left[ - \frac{1}{2} , \frac{1}{2} \right)^{d-1} \,\, \middle| \,\, |x_1 \cdots x_{d-1}| < \frac{c}{|x_d|} \right\} . \end{multline*}
If $|x_d|< c 2^{d-1}$, then the cross section has Lebesgue measure 1. Otherwise, using the inductive hypothesis, the Lebesgue measure of the cross section is
\[ O \left( \frac{c}{|x_d|} \log^{d-2} \frac{|x_d|}{c} \right) = O \left( \frac{c}{|x_d|} \log^{d-2} \frac{1}{c} \right) . \]
Applying Fubini's theorem we thus obtain
\[ \lambda (U_c) = c 2^d + O \left( c \log^{d-2} \frac{1}{c} \int_{\left( - \frac{1}{2} , - c2^{d-1} \right) \cup \left( c 2^{d-1} , \frac{1}{2} \right)} \frac{1}{|x_d|} \, \mathrm{d} x_d \right) = O \left( c \log^{d-1} \frac{1}{c} \right) . \]

Let $g : A_h \to \left[ - \frac{1}{2} , \frac{1}{2} \right)^d$ be defined as
\[ g(m) = \left( m \alpha_1 , \dots , m \alpha_d \right) \pmod{1} . \]
Note that $g$ is injective because of the irrationality of $\alpha_1, \dots, \alpha_d$, and $g(A_h) \subset U_{h L_M}$. It is easy to see that there exists a partition of $\left[ - \frac{1}{2} , \frac{1}{2} \right)^d$ into congruent axis parallel cubes with common side length in the open interval $\left( \frac{1}{2}L_M^{\frac{1}{d}}, L_M^{\frac{1}{d}} \right)$. Let $\mathcal{C}$ denote the family of cubes in such a partition. Every cube in $\mathcal{C}$ contains at most one point of $g(A_h)$. Indeed, otherwise there would exist $1 \le m < m' \le M$ such that $\left\| (m'-m) \alpha_i \right\| < L_M^{\frac{1}{d}}$ for every $1 \le i \le d$, and so
\[ \left\| (m'-m) \alpha_1 \right\| \cdots \left\| (m'-m) \alpha_d \right\| < L_M , \]
contradicting \eqref{lower}. Therefore the pigeonhole principle implies that
\begin{equation}\label{cardinality}
|A_h| \le \left| \left\{ C \in \mathcal{C} \,\, \middle| \,\, C \cap U_{hL_M} \neq \emptyset \right\} \right| .
\end{equation}

For an arbitrary $x \in U_{hL_M}$ consider the product
\[ \left( |x_1| + L_M^{\frac{1}{d}} \right) \cdots \left( |x_d| + L_M^{\frac{1}{d}} \right) . \]
When expanding this product let us estimate one of the terms as $|x_1 \cdots x_d| \le h L_M$, and all the other terms by simply using $|x_i| \le \frac{1}{2}$. This way we get
\begin{equation}\label{cube}
\left( |x_1| + L_M^{\frac{1}{d}} \right) \cdots \left( |x_d| + L_M^{\frac{1}{d}} \right) = O \left( h L_M + L_M^{\frac{1}{d}} \right)
\end{equation}
for any $x \in U_{h L_M}$ with an implied constant depending only on $d$. The estimate \eqref{cube} shows that
\begin{equation}\label{union}
\bigcup \left\{ C \in \mathcal{C} \,\, \middle| \,\, C \cap U_{hL_M} \neq \emptyset \right\} \subseteq U_c
\end{equation}
for some $c = O \left( h L_M + L_M^{\frac{1}{d}} \right)$. Comparing the Lebesgue measures of the sets in \eqref{union}, and using \eqref{cardinality} we get
\begin{equation}\label{bound}
|A_h| = O \left( \frac{\lambda (U_c)}{L_M} \right) = O \left( \left( h + L_M^{\frac{1}{d}-1} \right) \log^{d-1} \frac{1}{L_M} \right)
\end{equation}
with an implied constant depending only on $d$.

For every integer $k \ge 0$ let
\[ B_k = \left\{ 1 \le m \le M \,\, \middle| \,\, 2^k L_M \le \left\| m \alpha_1 \right\| \cdots \left\| m \alpha_d \right\| < 2^{k+1} L_M \right\} \subseteq A_{2^{k+1}} . \]
Note that if $k > \log_2 \frac{1}{L_M}$, then $B_k = \emptyset$. Therefore \eqref{bound} implies
\[ \begin{split} \sum_{m=1}^M \frac{1}{\left\| m \alpha_1 \right\| \cdots \left\| m \alpha_d \right\|} &\le \sum_{0 \le k \le \log_2 \frac{1}{L_M}} \frac{1}{2^k L_M} \left| A_{2^{k+1}} \right| \\ &= \left\{ \substack{ \displaystyle{ O \left( \frac{1}{L_M} \log \frac{1}{L_M} \right) \,\, \mathrm{if} \,\, d=1,} \\ \displaystyle{ O \left( \frac{1}{L_M^{2-\frac{1}{d}}} \log^{d-1} \frac{1}{L_M} \right) \,\, \mathrm{if} \,\, d \ge 2.}} \right. \end{split} \]

Finally, for arbitrary algebraic reals $\alpha_1, \dots , \alpha_d$ such that $1, \alpha_1, \dots , \alpha_d$ are linearly independent over $\mathbb{Q}$, Schmidt's theorem \eqref{Schmidt2} implies that $\frac{1}{L_M} = O \left( M^{1 + \varepsilon} \right)$ for any $\varepsilon >0$, and hence
\[ \sum_{m=1}^M \frac{1}{\left\| m \alpha_1 \right\| \cdots \left\| m \alpha_d \right\|} = O \left( M^{2-\frac{1}{d} + \varepsilon} \right) \]
for any $\varepsilon >0$.
\begin{flushright}
$\square$
\end{flushright}

\noindent\textbf{Proof of Theorem \ref{theorem6}:}

\noindent\textbf{(i)} We shall in fact prove that for any $1 \le T_1 < T_2$ and $\varepsilon >0$ we have
\begin{equation}\label{average}
\frac{1}{T_2-T_1} \int_{T_1}^{T_2} \left( \left| tC \cap \mathbb{Z}^d \right| - p(t) \right) \, \mathrm{d}t = O \left( 1+ \left( T_2 - T_1 \right)^{1 - \gamma_{d-1}- \varepsilon} \right) .
\end{equation}
Proposition \ref{proposition4} yields
\begin{multline}\label{bound2}
\frac{1}{T_2-T_1} \int_{T_1}^{T_2} \left( \left| tC \cap \mathbb{Z}^d \right| - p(t) \right) \, \mathrm{d}t \\ = \frac{1}{T_2 - T_1} \int_{T_1}^{T_2} E_N (t) \, \mathrm{d}t + O \left( 1 + T_2^{d-1+\varepsilon} \sqrt{\frac{\log N}{N}}  \right)
\end{multline}
for any integer $N>1$, where $E_N (t)$ is as in Definition \ref{definition3}. To estimate the average of $E_N (t)$ note that for any integer $m_j \neq 0$ we have
\begin{equation}\label{EN1}
\left| \frac{1}{T_2 - T_1} \int_{T_1}^{T_2} e^{-2 \pi i m_j a_j t} \, \mathrm{d}t \right| \le \min \left( 1 , \frac{1}{(T_2 - T_1) \pi |m_j| a_j} \right) .
\end{equation}
Indeed, using the triangle inequality we get that the left hand side of \eqref{EN1} is at most $1$. On the other hand, by explicitly evaluating the integral we get
\[ \left| \frac{1}{T_2 - T_1} \int_{T_1}^{T_2} e^{-2 \pi i m_j a_j t} \, \mathrm{d}t \right| = \frac{|e^{-2 \pi i m_j a_j T_2} -e^{-2 \pi i m_j a_j T_1}|}{(T_2-T_1)2 \pi |m_j| a_j} , \]
where the numerator is clearly at most $2$. Elementary calculation shows that for any $c \in \mathbb{R} \backslash \mathbb{Z}$ and any integer $M \ge 0$ we have the general estimate
\begin{equation}\label{general}
\begin{split} \left| \sum_{m=-M}^M \frac{1}{c-m} \right| &\le \frac{1}{|c|} + \sum_{m=1}^M \frac{2 |c|}{|c^2 - m^2|} \\ &\le \frac{1}{|c|} + \sum_{1 \le m < |c|-1} \frac{2|c|}{(|c|-m)(|c|+m)} + \frac{3}{\left\| c \right\|} + \sum_{m>|c|+1} \frac{2 |c|}{m^2 - c^2} \\ &= O \left( \frac{\log (|c|+1)}{\left\| c \right\|} \right) . \end{split}
\end{equation}
Applying \eqref{general} with $c=m_j \frac{a_j}{a_k}$ and $M=M_k$ for every $k \neq j$, and using \eqref{EN1} we obtain
\begin{equation}\label{EN2}
\frac{1}{T_2 - T_1} \int_{T_1}^{T_2} E_N (t) \, \mathrm{d}t = O \left( \sum_{j=1}^d \sum_{m=1}^{\infty} \frac{\log^{d-1} (m+1)}{m \prod_{k \neq j} \left\| m \frac{a_j}{a_k} \right\|} \min \left( 1 , \frac{1}{(T_2 - T_1) m} \right) \right) .
\end{equation}

Let us first estimate the terms $1 \le m \le \frac{1}{T_2 - T_1}$. Using Definition \ref{definition4} of $\gamma_d$ we get that for any integer $\ell \ge 0$ we have
\[ \sum_{2^{\ell} \le m < 2^{\ell +1}} \frac{\log^{d-1} (m+1)}{m \prod_{k \neq j} \left\| m \frac{a_j}{a_k} \right\|} = O \left( \frac{(\ell +1)^{d-1}}{2^{\ell}} \left( 2^{\ell +1} \right)^{\gamma_{d-1} + \varepsilon} \right) . \]
Summing over $0 \le \ell \le \log_2 \frac{1}{T_2 - T_1}$ we obtain
\begin{equation}\label{est1}
\sum_{1 \le m \le \frac{1}{T_2 - T_1}} \frac{\log^{d-1} (m+1)}{m \prod_{k \neq j} \left\| m \frac{a_j}{a_k} \right\|} = O \left( (T_2 - T_1)^{1 - \gamma_{d-1} - \varepsilon} \right) .
\end{equation}
To estimate the terms $m \ge \frac{1}{T_2 - T_1}$ let again $\ell \ge 0$ be an integer and consider
\[ \sum_{2^{\ell} \le m < 2^{\ell +1}} \frac{\log^{d-1} (m+1)}{m^2 \prod_{k \neq j} \left\| m \frac{a_j}{a_k} \right\| (T_2 - T_1)} = O \left( \frac{(\ell + 1)^{d-1}}{2^{2 \ell} (T_2 - T_1)} \left( 2^{\ell +1} \right)^{\gamma_{d-1} + \varepsilon}\right) . \]
Using the fact $\gamma_{d-1} <2$ from Theorem \ref{theorem5} we can sum over every $\ell \ge \log_2 \frac{1}{T_2 - T_1}-1$ to obtain
\begin{equation}\label{est2}
\sum_{m \ge \frac{1}{T_2 - T_1}} \frac{\log^{d-1} (m+1)}{m^2 \prod_{k \neq j} \left\| m \frac{a_j}{a_k} \right\| (T_2 - T_1)} = O \left( (T_2 - T_1)^{1 - \gamma_{d-1} - \varepsilon} \right) .
\end{equation}
Thus \eqref{est1} and \eqref{est2} imply that \eqref{EN2} simplifies as
\begin{equation}\label{EN3}
\frac{1}{T_2 - T_1} \int_{T_1}^{T_2} E_N (t) \, \mathrm{d}t = O \left( (T_2 - T_1)^{1-\gamma_{d-1} - \varepsilon} \right) .
\end{equation}
Using \eqref{EN3} in \eqref{bound2}, and letting $N \to \infty$ we obtain \eqref{average}, as claimed.

\vspace{5mm}

\noindent\textbf{(ii)} The main idea is to use the fact that $\left| tC \cap \mathbb{Z}^d \right|$ is a monotone nondecreasing function of the real variable $t>1$. Fix a real number $t>1$. Since $p$ is a polynomial of degree $d$, there exists a constant $K_1>1$ such that
\[ \left| p(t+h) - p(t) \right| \le K_1 t^{d-1} |h| \]
for any $-1<h<1$. Let $\Delta (t) = \left| tC \cap \mathbb{Z}^d \right| -p(t)$. The trivial bound \eqref{eq2} gives that $\left| \Delta (t) \right| \le K_2 t^{d-1}$ for some constant $K_2>1$. Let $K = \max \left\{ K_1 , K_2 \right\}$.

If $\Delta (t) >0$, then for any $u \in \left[ t, t+ \frac{\Delta (t)}{2K t^{d-1}} \right]$ we have
\[ \left| uC \cap \mathbb{Z}^d \right| - p(u) \ge \left| tC \cap \mathbb{Z}^d \right| - p(t) - (p(u) - p(t)) \ge \frac{\Delta (t)}{2} . \]
Applying \eqref{average} from (i) to the interval $[T_1 , T_2] = \left[ t, t+ \frac{\Delta (t)}{2K t^{d-1}} \right]$ we obtain
\begin{equation}\label{delta}
\Delta (t) = O \left( \left( \frac{| \Delta (t)|}{t^{d-1}} \right)^{1-\gamma_{d-1} - \varepsilon} \right) .
\end{equation}
Similarly, if $\Delta (t) < 0$, then for any $u \in \left[t- \frac{\Delta (t)}{2Kt^{d-1}} ,t \right]$ we have
\[ \left| uC \cap \mathbb{Z}^d \right| -p(u) \le \left| tC \cap \mathbb{Z}^d \right| -p(t) + (p(t) -p(u)) \le \frac{\Delta (t)}{2} . \]
Applying \eqref{average} from (i) to the interval $[T_1 , T_2] = \left[ t- \frac{\Delta (t)}{2K t^{d-1}} ,t \right]$ we obtain that \eqref{delta} holds in the case $\Delta (t) <0$ as well. Rearranging \eqref{delta} we get
\[ \Delta (t) = O \left( t^{\frac{\gamma_{d-1}-1}{\gamma_{d-1}} (d-1) + \varepsilon} \right) \]
for any $\varepsilon >0$, as claimed.
\begin{flushright}
$\square$
\end{flushright}

\noindent\textbf{Proof of Theorem \ref{theorem7}:} Proposition \ref{proposition1} and Definition \ref{definition4} yield
\[ \left| tS \cap \mathbb{Z}^d \right| - q(t) = \frac{1}{2^d} \sum_{I \subseteq [d]} \left( \left| t C_I \cap \mathbb{Z}^d \right| - p_I (t) \right) , \]
where $p_I = p_{\left( a_i \,\, \middle| \,\, i \in I \right) }$. Since the terms with $|I| \le 1$ can be estimated easily, we can reduce Theorem \ref{theorem7} to Theorem \ref{theorem6} in dimensions $2, 3,  \dots , d$. It is easy to see from Definition \ref{definition4} that $\gamma_1 \le \gamma_2 \le \cdots \le \gamma_{d-1}$, and so
\[ \frac{\gamma_{|I|-1}-1}{\gamma_{|I|-1}} \left( |I| -1 \right) \le \frac{\gamma_{d-1}-1}{\gamma_{d-1}} (d-1) \]
for any $2 \le |I| \le d$.
\begin{flushright}
$\square$
\end{flushright}

\end{document}